\documentclass[12pt,leqno]{article}

\usepackage{amssymb}
\usepackage{amsmath}
\usepackage{graphicx,psfrag,color}
\usepackage{url}
\usepackage{pstool}

\newtheorem{thm}{Theorem}
\newtheorem{lem}{Lemma}

\newtheorem{Def}{Definition}
\newtheorem{qes}{Question}

\newcommand{\Gg}{\mbox{$\mathfrak G$}}
\newcommand{\Mm}{\mbox{$\mathfrak M$}}
\newcommand{\Nn}{\mbox{$\mathfrak N$}}

\newcommand{\Aa}{\mbox{$\mathfrak A$}}

\newcommand{\Ll}{\mathcal{L}}
\newcommand{\Th}{\mbox{\sf Th}}
\newcommand{\Type}{\mbox{\sf Type}}

\newcommand{\ty}{\mbox{\sf ty}}
\newcommand{\Types}{\mbox{\sf Types}}

\newcommand{\inv}[1]{{#1}^{\smallsmile}}
\newcommand{\id}[1]{\mbox{\sf id}({#1})}
\newcommand{\RA}[1]{\mbox{$\mathfrak{R}\mathfrak{a}(#1)$}}
\newcommand{\nequiv}{{\overset{\small n}{\equiv}}}
\newcommand{\tequiv}{{\overset{\small 2}{\equiv}}}
\newcommand{\R}{\mbox{\sf R}}

\newcommand{\Qed}{\hfill\mbox{$\Box$}\bigskip}

\newcommand{\comment}[1]{}

\begin{document}

\title{Two-variable logic has weak, but not strong, Beth definability}
\author{Andr\'eka, H.\ and N\'emeti, I.}%
\date{November 2020}
\maketitle

\begin{abstract} We prove that the two-variable fragment of
first-order logic has the weak Beth definability property. This
makes the two-variable fragment a natural logic separating the weak
and the strong Beth properties since it does not have the strong
Beth definability property.
\end{abstract}

\begin{center}
\emph{We dedicate this paper to Harvey Friedman in respect for his
work}
\end{center}

\section{Introduction}\label{intro-sec}

One of the many expressibility properties of first-order logic with
equality FO is the Beth definability property BDP. It states that if
a relation can be specified by some extra means then it can be
specified explicitly without using the extra means. In more detail,
if $\Th$ is a first-order logic theory in a language $\Ll$ and
$\Sigma$ is another first-order logic theory in the language $\Ll$
expanded with an extra relation symbol $\R$ such that in each model
of $\Th$ there is at most one relation $R$ satisfying $\Sigma$, then
this unique relation can be defined in the original language $\Ll$
without using the extra relation symbol, i.e., there is a formula
$\varphi$ in $\Ll$ such that $\Th\cup\Sigma\models
\forall\overline{x}(\R(\overline{x})\leftrightarrow\varphi)$. In
this context, $\Sigma$ is called the implicit definition and
$\varphi$ is called the explicit definition of $\R$.

Investigating BDP for fragments of FO means showing that if all the
formulas of the theory and the implicit definition belong to the
fragment, then the explicit definition, too, belongs to it. Thus,
having BDP or not shows a kind of ``integrity" of the fragment, and
a kind of ``complexity-property" of FO itself.
For example, the guarded fragment GF of FO has BDP \cite{HMO99}.
Thus, if the theory $\Th$ and the implicit definition $\Sigma$
consist of guarded formulas, then the explicit definition can be
chosen to be guarded, too. We note that in the strict sense, GF does
not have Craig interpolation property, yet it is worth deciding when
the interpolant belonging to guarded $\Th$ and $\Sigma$ can be
chosen to be guarded itself \cite{JW20}. For work of the similar
kind, see e.g., \cite{BBtC}.

It is known that $n$-variable fragments FO$n$ of FO do not have BDP,
for all finite $n\ge 2$, see \cite{ACMNS}. 
This means that if the theory $\Th$ and the implicit definition
$\Sigma$ all use only $n$ variables, the explicit definition
$\varphi$ may need more than $n$ variables. That the BDP fails for
FO$2$ is kind of surprising, because FO$2$ usually behaves ``better"
than FO$n$ for $n\ge 3$. For example, FO$2$ is decidable while
FO$n$, for $n\ge 3$ is not.

The weak Beth definability property wBDP was introduced by Harvey
Friedman \cite{HFrie}. The definition of wBDP is the same as that of
BDP except that only those implicit definitions have to be made
explicit which also have the existence property, not only the
uniqueness property. Let us call these strong implicit definitions.
Thus we require that $\R$ has an explicit definition only when in
each model of $\Th$ there is exactly one relation $R$ satisfying
$\Sigma$, as opposed to having at most one such relation. (For
formal definition see Definition~\ref{wbeth-def} at the beginning of
section~\ref{wbeth-sec}.) Thus, BDP implies wBDP, since wBDP
requires fewer definitions to be equivalent to an explicit one. In
mathematical practice, one almost always requires both existence and
uniqueness for an implicitly defined object. For this reason, wBDP
sometimes is considered to be a more natural definability property
than BDP itself (see, e.g., \cite[p.129]{MS85}). We note that wBDP
was intensely investigated in abstract model theory in connection
with logics stronger than FO, see, e.g., \cite{BF, MSS}.

It is known that FO$n$ does not have wBDP either, whenever $n\ge 3$
(\cite{SS} for $n=3$, \cite{Hodk} for $n\ge 5$, and \cite{ANwBeth14}
for $n\ge 3$). Proving failure of wBDP amounted to find also strong
implicit definitions in FO$n$ that could not be made explicit. It
remained open whether FO$2$ had wBDP or not.

In this paper we prove that FO$2$ does have wBDP, that is, in FO$2$
all strong implicit definitions can be made explicit. This restores
two-variable logic's image that it behaves better than $n$-variable
logics for $n\ge 3$. This theorem may also point to wBDP being a
more natural property than BDP.

As far as we know, the present paper contains the first proof for a
logic to have the wBDP not via showing that it has the stronger BDP.
So, the difference between BDP and wBDP was not tangible so far in
the sense that there was no example for a logic that distinguished
the two properties.
FO augmented with the quantifier ``there exists uncountable many"
$L(Q)$ was a good candidate for such a distinguishing logic, since
it does not have BDP \cite{HFrie} and it is consistent with set
theory that it has wBDP \cite{MS85}. However, it is still an open
problem whether wBDP can be proved for $L(Q)$ in set theory or not.

It is also satisfying that the distinguishing logic FO$2$ is a
well-investigated, natural logic. Luckily, we did not have to
construct a logic to show that wBDP and BDP are distinct properties,
a well-known logic turned out to do the job for us. Two-variable
logic and its extensions are quite popular in computer science and
in modal logic.

Our proof hinges on the fact that FO$2$ has a special property that
FO$n$ with $n\ge 3$ do not have. Namely, each model of FO$2$ is
FO$2$-equivalent with a 
model in which elements of the same FO$2$-type also  have the same
automorphism-type (Theorem~\ref{trans-thm},
Theorem~\ref{trans2-thm}). FO$3$ does not have this property
(Theorem~\ref{notrans-thm}).

The structure of the paper is as follows. In section~\ref{trans-sec}
we define the above kind of models that we will call transitive and
we prove the key property of FO$2$ about the abundance of transitive
models. In section~\ref{wbeth-sec} we prove that FO$2$ has wBDP by
relying on these transitive models (Theorem~\ref{wbeth-thm}). We
close the paper with a remark about some connections with algebra
and the literature.

\section{Transitive models of FO$2$}\label{trans-sec}

We use the notation of \cite{CK}, if not stated otherwise. By FO we
mean first-order logic with equality, but we do not allow function
or constant symbols.

Let $n$ be a finite number. By FO$n$ we mean the fragment of FO that
uses only the first $n$ variables. Strictly speaking, the fragment
FO$n$ of FO is defined by taking for all languages $\Ll$ all the
models of $\Ll$ but restricting the set of formulas of FO to those
that contain the first $n$ variables only. Thus, relation symbols of
arbitrarily high rank can be allowed in FO$n$. For simplicity, in
this paper in FO$n$ we will allow only languages with relation
symbols of rank at most $n$. We go further, we allow relation
symbols of rank $n$ only. These are not important restrictions.
Usually, we do not indicate the language, but we will always work
with similar models, i.e., models having the same language, if not
stated otherwise.

\comment{In FO, there is a good notion for a homogeneous model.
However, for FO$n$, the natural $n$-variable analogues of several
equivalent properties of homogeneous models become non-equivalent,
this is one of the reasons we will not call our models
``homogeneous" but rather transitive ones.}

We say that $\Mm$ and $\Nn$ are \emph{$n$-equivalent}, in symbols
$\Mm\nequiv\Nn$, when the same n-variable formulas are true in $\Mm$
and in $\Nn$.
In this paper, we concentrate on $n=2$. By the \emph{type} of an
element in a model we understand the set of FO$2$-formulas true of
it. We say that $\Mm$ is \emph{transitive} if whenever $a$ and $b$
are elements having the same type in it, there is an automorphism of
$\Mm$ taking $a$ to $b$. In concise form:
\[\mbox{$(\Mm,a)\tequiv (\Mm,b)$ implies
$(\Mm,a)\cong (\Mm,b).$}\]
A typical transitive model is the set of integers with the successor
relation. A typical non-transitive model is a connected graph on a
set with more than 3 elements in which there are nodes of distinct
degree.

Let us call a model \emph{binary} if all its basic relations are of
rank 2. We prove in this section that each binary model is
2-equivalent to a transitive model (Theorem~\ref{trans2-thm}). The
idea of the proof is that we construct a 2-equivalent version of any
binary model via replacing each binary relation in it with a
suitable set of different successor relations.

A stronger version of the above will be proved in this section.
We call a model $\Mm$ \emph{$2$-homogeneous} if whenever $a$ and $b$
have the same type in it, to any element $c\in M$ there is $d\in M$
such that $(\Mm,a,c)\tequiv(\Mm,b,d)$. This is a straightforward
analogue of the definition of $\alpha$-homogeneity in \cite{CK}
where $\alpha$ is any ordinal. We note that transitivity implies
$2$-homo\-gen\-eity, but not the other way round. Further, a finite
model is always 2-homogeneous (see the proof of
Theorem~\ref{trans2-thm}).

Finally, we need the notion of $2$-partial isomorphism.
The notion of $n$-partial isomorphism was defined in
\cite[p.259]{Bent} as a natural restriction of the usual notion of
partial isomorphisms between models of FO (see \cite{CK}). We recall
the definition of 2-partial isomorphism in detail because we will
rely on it.

\begin{Def}[2-partial isomorphism]\label{piso-def}
The set $I$ is a \emph{2-partial isomorphism} between models
$\Mm, \Nn$ if (i)-(iv) below hold:
\begin{description}
\item{(i)}
$I$ relates elements as well as pairs of $M$ and $N$, i.e., it is a
subset of $(M\times N)\cup (M^2\times N^2)$,
\item{(ii)} local isomorphism property:\\
related pairs of $I$ are isomorphisms between $\Mm$ and $\Nn$
restricted to the first and second parts of the pair, respectively,
\item{(iii)} restriction property:\\ if $\langle
(a,a'),(b,b')\rangle\in I$ then $\langle a,b\rangle\in I$ and
$\langle a',b'\rangle\in I$,  and
\item{(iv)}
back-and-forth property:\\
$\forall a\in M\exists b\in N\langle a,b\rangle\in I$ and vice
versa, $\forall b\in N\exists a\in M\langle a,b\rangle\in
I$,\\
$\forall\langle a,b\rangle\in I\forall a'\in M\exists b'\in
N\langle(a,a'),(b,b')\rangle\in I$, and vice
versa\\
$\forall\langle a,b\rangle\in I\forall b'\in N\exists a'\in
M\langle(a,a'),(b,b')\rangle\in I$. \Qed\end{description}
\end{Def}
Instead of 2-partial isomorphism we will simply say
\emph{2-isomorphism}. It is known that if two pairs are related by a
2-isomorphism, then the same FO$2$-formulas are true of them, this
is straightforward to show by induction. Thus, if there is a
2-isomorphism between two models, then they are 2-equivalent.

Being 2-isomorphic is stronger than being 2-equivalent. We are going
to prove that a binary model is 2-isomorphic to a transitive model
if and only if it is 2-homogeneous.
This stronger theorem will be used in the proof of weak Beth
definability property for 2-variable logic FO$2$
(Theorem~\ref{wbeth-thm}).

\begin{thm}
\label{trans-thm} A binary model $\Mm$ is 2-isomorphic to a
transitive model $\Nn$ if and only if $\Mm$ is 2-homogeneous.
\end{thm}

\noindent{\bf Proof.} We prove sufficiency of 2-homogeneity first.
Let $\Mm$ be any 2-homogeneous binary model. We are going to define
another model $\Nn$ and a 2-isomorphism $I$ between them. $\Nn$ will
be finite when $\Mm$ is so. Then we show that $\Nn$ is transitive.

Some notation and terminology. The variables of FO$2$ will be
denoted by $x,y$. By the \emph{2-type} of $a,b\in \Mm$ we understand
the set of FO$2$-formulas $\rho(x,y)$ that are true for $a,b$ in
$\Mm$. Formally,
\[ \Type(a,b,\Mm) = \{\rho(x,y)\in FO2 : \Mm\models\rho[a,b]\} .\]
Now we define
\begin{description}
\item{} $[a] = \{ b\in M : \Type(b,b,\Mm)=\Type(a,a,\Mm)\}$,
\item{} $[a,b] = \Type(a,b,\Mm)$,
\item{} $\Types([a],[b]) = \{ [p,q] : p\in[a],\ q\in[b]\}$.
\end{description}
We call the elements of $\{[a] : a\in M\}$ \emph{1-Types}. Note that
$[a]$ is a subset of the model, while $[p,q]$ is a set of
FO$2$-formulas.

The type $[p,q]$ determines $[q,p]$, we call the latter the
\emph{converse type} $\inv{[p,q]}$ of $[p,q]$. Taking the converse
is a bijection between $\Types([a],[b])$ and $\Types([b],[a])$. We
call a type $[p,q]$ \emph{symmetric} when $[p,q]=[q,p]$, otherwise
we call it \emph{asymmetric}. There is a unique element of
$\Types([a],[a])$ which contains $x=y$, we call this the
\emph{identity type on $[a]$}, it is denoted by $\id{a}$ and it is
symmetric. A \emph{non-identity} type is a type which is not an
identity type.

We begin the definition of $\Nn$.
A group $\Gg=\langle G,+\rangle$ will be used in the definition of
$\Nn$. That is, $G$ is the universe of the group, and $+$ is its
group-operation. The identity element (or zero-element) of $+$ is
denoted by $0$, and the inverse of an element $a$ is denoted by
$-a$. We call a group \emph{asymmetric} when its zero-element is its
only element of order 2, i.e., when $x+x=0$ implies $x=0$ in it.
Let $\Gg$ be any commutative asymmetric group of size at least twice
that of $\Types([a],[b])$ for any $a,b\in M$, i.e.,
\[\tag{g1}\label{g1} |G| \ge 2\cdot|\Types([a],[b])|\quad\mbox{ for all $a,b\in
M$} .\]  There is such a $\Gg$ because for each odd number $n$ the
group $\mathbb{Z}_n$ of integers smaller than $n$ and with addition
modulo $n$ is asymmetric, and then one can construct an asymmetric
group of any infinite size by taking an elementary submodel of a
sufficiently big ultraproduct of the $\mathbb{Z}_n$ with odd $n$.
Equivalently, use the upward L\"owenheim-Skolem-Tarski Theorem
\cite[Cor.2.1.5, Cor.2.1.6]{CK}.

The universe $N$ of $\Nn$ is defined as
\[ N = \{ ([a],g) : a\in M, g\in G\} .\]
Clearly, if $M$ is finite, then $G$ can be chosen to be finite, and
then $N$ is finite.

For the definition of the relations of $\Nn$, let
$\lambda=\langle\lambda_{[a],[b]} : a,b\in M\rangle$ be a system of
functions mapping $G$ to the types of $\Mm$ that satisfies the
following conditions for all $a,b\in M$. We will write
$\lambda_{a,b}$ in place of $\lambda_{[a],[b]}$, for easier
readability.
\begin{description}
\item{} $\lambda_{a,b}:G\to\Types([a],[b])$ is surjective,
\item{} $\lambda_{a,b}(g) = \inv{\lambda_{b,a}(-g)}\quad\mbox{and}\quad
\lambda_{a,a}(0)=\id{a}$.
\end{description}
There is such a system $\lambda$ of functions, because of the
following. When $[a]\ne[b]$, take any surjective
$\lambda_{b,a}:G\to\Types([b],[a])$, there is such since
$|G|\ge|\Types([b],[a])|$ by (g1), then define
$\lambda_{a,b}(g)=\inv{\lambda_{b,a}(-g)}$,  this is also surjective
because $\Types([a],[b])=\{ \inv{-t} : t\in\Types([b],[a])\}$.

Assume now $[a]=[b]$. Since $G$ is asymmetric, there is a set
$P\subseteq G$ such that  $P$, $-P=\{ -g : g\in P\}$ together with
$\{ 0\}$ form a partition of $G$. Let $S$ denote the set of all
non-identity symmetric types in $T=\Types([a],[a])$, then there is a
set $A\subseteq T$ such that $A$, $\inv{A}=\{\inv{t} : t\in A\}$,
and $S$ together with $\{\id{a}\}$ form a partition of $T$. Now,
$|P|\ge|S\cup A|$ by (g1), so there is a surjective function
$L:P\to(S\cup A)$. Define now $\lambda_{a,a}(g)=L(g)$ for $g\in P$,
$\lambda_{a,a}(g)=\inv{L(-g)}$ for $g\in -P$ and
$\lambda_{a,a}(0)=\id{a}$. This function satisfies the required
conditions (because both functions of taking inverse in the group
and taking converse in the types are their own inverses).

Now we define, for all $a,b\in M$ and $g,h\in G$
\[ \ty(([a],g),([b],h)) = \lambda_{a,b}(h-g), \]
where $h-g$ denotes $h+-g$ as usual in group theory. Then $\ty$ maps
$N\times N$ to the types of $\Mm$. Let $\R$ be an arbitrary binary
relation symbol in the language of $\Mm$. The binary relation
$R^{\Nn}$ belonging to  $\R$ in $\Nn$ is defined as
\[ R^{\Nn} = \{ (p,q)\in N\times N : \R(x,y)\in\ty(p,q)\} .\]
By this, the model $\Nn$ has been defined.

Next, we exhibit a 2-isomorphism between $\Mm$ and $\Nn$. We define
$I\subseteq (M\times N)\cup (M^2\times N^2)$ by requiring for all
$a,b\in M$ and $p,q\in N$ that

\begin{description}
\item{}
$\langle a,p\rangle\in I\quad\mbox{ iff }\quad p=([a],g)$ for some
$g$,
\item{}
$\langle (a,b),(p,q)\rangle\in I\quad\mbox{ iff }\quad
[a,b]=\ty(p,q)$.
\end{description}

\noindent We now show that $I$ is a 2-isomorphism between $\Mm$ and
$\Nn$. From the conditions defining a 2-isomorphism, $I$ clearly
satisfies (i) by its very definition.

Next we show that the restriction property (iii) holds for $I$.
Assume that $\langle (a,b),(p,q)\rangle\in I$. This means that
$[a,b]=\ty(p,q)$. Assume that $p=([c],g)$ and $q=([d],h)$. From the
definition of $\ty(p,q)$ it is clear that
$\ty(p,q)\in\Types([c],[d])$. Thus, $[a,b]=[r,s]$ for some $r\in[c]$
and $s\in[d]$. But this implies that $[a]=[r]=[c]$ and
$[b]=[s]=[d]$, hence $\langle a,p\rangle\in I$ and $\langle
b,q\rangle\in I$.

We show that $I$ satisfies local isomorphism property (ii).  Let
$\R$ be an arbitrary relation symbol in the language of $\Mm$.
Assume that $\langle a,p\rangle\in I$. Then $p=([a],g)$ for some $g$
by the definition of $I$, and we have to show that $R(p,p)$ holds in
$\Nn$ iff $R(a,a)$ holds in $\Mm$. By definition, $R(p,p)$ holds in
$\Nn$ iff
$\R(x,y)\in\ty(p,p)=\lambda_{a,a}(g-g)=\lambda_{a,a}(0)=\id{a}$, and
$\R(x,y)\in\id{a}$ iff $R(a,a)$ holds in $\Mm$, by the definition of
$\id{a}$.

Assume that $\langle (a,b),(p,q)\rangle\in I$. Then $[a,b]=\ty(p,q)$
by the definition of $I$. We have to show that each of $a=b,
R(a,a)$, $R(a,b)$, $R(b,a)$, and $R(b,b)$ holds in $\Mm$ iff the
same holds in $\Nn$ for $p,q$ in place of $a,b$. Assume that
$p=([c],g)$ and $q=([d],h)$, then $\ty(p,q)=\lambda_{c,d}(h-g)$ by
the definition of $\ty$. Now, $a=b$ iff $x=y\in
[a,b]=\ty(p,q)=\lambda_{c,d}(h-g)$, and $x=y\in\lambda_{c,d}(h-g)$
iff ($[c]=[d]$ and $h-g=0$) iff $p=q$. Similarly, $R(a,b)$ iff
$\R(x,y)\in[a,b]=\ty(p,q)$ iff $R(p,q)$, by the definition of $\Nn$.
For the next case, we want to show that
\[\tag{t1} \ty(p,q)=\inv{\ty(q,p)}.\] Indeed,
$\ty(p,q)=\lambda_{c,d}(h-g)=\inv{\lambda_{d,c}(g-h)}=\inv{\ty(q,p)}$,
by the second condition that $\lambda$ has to satisfy. By this, (t1)
is proved.

Now, $R(b,a)$ iff
$\R(x,y)\in[b,a]=\inv{[a,b]}=\inv{\ty(p,q)}=\ty(q,p)$ iff $R(q,p)$.
Finally, $R(a,a)$ iff $R(p,p)$ and $R(b,b)$ iff $R(q,q)$ hold by the
first case of (ii), since we have already shown the restriction
property (iii). Thus, $\langle (a,b),(p,q)\rangle$ indeed specifies
a partial isomorphism between $\Mm$  restricted to $\{ a,b\}$  and
$\Nn$ restricted to $\{ p,q\}$.

We check the back-and-forth property (iv) for $I$. To check the
first part, notice that to any $a\in M$ there is at least one
$p=([a],g)$ in $N$, because $G$ is nonempty. The second part is
clear, since no 1-Type $[a]$  is empty. To check the third and
fourth parts of the back-and-forth property, assume that $\langle
a,p\rangle\in I$ with $p=([a],g)$. Let $b\in M$ be arbitrary. We
have to find a $q\in N$ such that $[a,b]=\ty(p,q)$. By surjectivity
of $\lambda_{a,b}$, there is $f\in G$ with $[a,b]=\lambda_{a,b}(f)$.
Let $q=(b,f+g)$. Then
$\ty(p,q)=\lambda_{a,b}(f+g-g)=\lambda_{a,b}(f)=[a,b]$ and we are
done with the third part.
Let now $q=([b],h)\in N$ be arbitrary. We have to find $c\in M$ such
that $[a,c]=\ty(p,q)$. Now, $\ty(p,q)\in\Types([a],[b])$ which means
that there are $a'\in[a]$ and $b'\in[b]$ such that
$\ty(p,q)=\Type(a',b',\Mm)$. By $a'\in[a]$ we have
$\Type(a',a',\Mm)=\Type(a,a,\Mm)$, so by 2-homogeneity of $\Mm$
there is $c\in M$ such that $\Type(a,c,\Mm)=\Type(a',b',\Mm)$ and
this shows that $[a,c]=\ty(p,q)$. So, (iv) holds for $I$.
 We have seen that $I$ is a 2-isomorphism between $\Mm$ and
$\Nn$.

We show that $\Nn$ is transitive. We have to show that if $p,q\in N$
are of the same type in $\Nn$, then there is an automorphism of
$\Nn$ that takes $p$ to $q$. Assume that $p=([a],g)$ and
$q=([b],h)$. Then $\langle a,p\rangle\in I$ and $\langle
b,q\rangle\in I$, by the definition of $I$.  Since $I$ is a
2-isomorphism between $\Mm$ and $\Nn$, we get that the same formulas
hold in $\Nn$ for $p$ as in $\Mm$ for $a$, and the same for $q$ and
$b$. Hence, $a$ and $b$ are of the same type in $\Mm$ (since $p$ and
$q$ are of the same type in $\Nn$).
It is easy to check that $a$ and $b$ are of the same type in $\Nn$
iff $[a]=[b]$. Let $k=-g+h$ and define $\alpha:N\to N$ by
\[ \alpha(([c],f)) = ([c],f+k)\quad\mbox{ for all $c\in M$ and $f\in G$}.\]
Now, $\alpha(p)=q$ by $[a]=[b]$ and $h=g+k$. We show that $\alpha$
is an automorphism of $\Nn$. First, $\alpha$ is a permutation of $N$
because $\Gg$ is a group. Let $\R$ be a binary relation symbol in
the language of $\Nn$, and let $r=([c],i)\in N$, $s=([d],j)\in N$.
Now, $R(r,s)$ holds in $\Nn$ iff
$\R(x,y)\in\ty(r,s)=\lambda_{c,d}(j-i)$. Similarly,
$R(\alpha(r),\alpha(s))$ holds in $\Nn$ iff
$\R(x,y)\in\ty(\alpha(r),\alpha(s))=\lambda_{c,d}(j+k-(i+k))$.
However, $j-i=j+k-(i+k)$ because $\Gg$ is commutative. Thus,
$\alpha$ is indeed an automorphism, and we are done with showing
that $\Nn$ is transitive.

We have seen that $\Nn$ is transitive, and this finishes the proof
of one direction of Theorem~\ref{trans-thm}.

The other direction of Theorem~\ref{trans-thm}, necessity of
2-homogeneity, follows from the facts that a transitive model is
always 2-homogeneous, and 2-isomorphisms preserve being
2-homogeneous. In more detail, assume that $I$ is a 2-isomorphism
between $\Mm$ and the transitive $\Nn$. We have to show that $\Mm$
is 2-homogeneous. Let $a,b,c\in M$ be such that $[a]=[b]$. By the
back-and-forth property in the definition of a 2-isomorphism, there
are $a',b',c'\in N$ such that $\langle a,a'\rangle$, $\langle
b,b'\rangle$, $\langle (a,c),(a',c')\rangle$ are all in $I$. Then
$\Type(a',a',\Nn)=\Type(b',b',\Nn)$ since $I$ is a 2-isomorphism and
$[a]=[b]$. Since $\Nn$ is 2-transitive, there is an automorphism of
$\Nn$ that takes $a'$ to $b'$. Let $d'$ be the image of $c'$ under
this automorphism. Then $\Type(a',c',\Nn)=\Type(b',d',\Nn)$ since
automorphisms preserve 2-types of elements. By the back-and-forth
property of $I$ again, there is $d\in M$ such that $\langle
(b,d),(b',d')\rangle\in I$. Then
$\Type(a,c,\Mm)=\Type(a',c',\Nn)=\Type(b',d',\Nn)=\Type(b,d,\Mm)$
and we are done.
 \Qed

Next we state a corollary of Theorem~\ref{trans-thm}.

\begin{thm}\label{trans2-thm} Each binary model is 2-equivalent to
a transitive model. Each finite binary model is 2-isomorphic to a
finite transitive model.
\end{thm}

\noindent{\bf Proof.} For the definition of an $\omega$-saturated
model see, e.g., \cite{CK}. First we prove
\[\tag{S} \mbox{$\Mm$ is $\omega$-saturated implies that $\Mm$ is
2-homogeneous.}\]
Assume that $\Mm$ is $\omega$-saturated. Let $a,b,c\in M$ be such
that $\Type(a,a,\Mm)=\Type(b,b,\Mm)$. Let $Y=\{ b\}\subseteq M$ and
let $\Gamma(x) = \{ \rho(x,b)\in FO$2$ : \rho(c,a)\mbox{ in }\Mm\}$.
Then $\Gamma(x)$ is a set of formulas in the language of
$\langle\Mm,b\rangle$. We show that it is consistent with the theory
of $\langle\Mm,b\rangle$. Let $\Delta$ be a finite subset of
$\Gamma(x)$, let $\delta(y)$ denote the formula $\exists
x\bigwedge\Delta[b/y]$ that we get from $\exists x\bigwedge\Delta$
by replacing $b$ everywhere with  $y$. Then $\Mm\models\delta(a)$ by
the definition of $\Gamma(x)$, and so $\Mm\models\delta(b)$ since
$a,b$ have the same 1-Type in $\Mm$. But $\Mm\models\delta(b)$ means
that $\langle\Mm,b\rangle\models\exists x\bigwedge\Delta$ that shows
that $\Delta$ is consistent with the theory of $\langle
\Mm,b\rangle$. Since $\Delta$ is an arbitrary finite subset of
$\Gamma(x)$, this means that $\Gamma(x)$ is consistent with the
theory of $\langle\Mm,b\rangle$. Since $\Mm$ is $\omega$-saturated,
then there is $d\in M$ for which $\Gamma(d)$ holds. This means that
$\Type(a,c,\Mm)=\Type(b,d,\Mm)$ and we are done with showing (S).

Now, Theorem~\ref{trans2-thm} follows from (S) and
Theorem~\ref{trans-thm} by using that each infinite model is
elementarily equivalent--hence 2-equivalent--with an
$\omega$-saturated one (see \cite[Lemma 5.1.4]{CK}), that each
finite model is $\omega$-saturated (see \cite[Prop.5.1.2]{CK}) and
that the model $\Nn$ in the proof of Theorem~\ref{trans-thm} can be
constructed to be finite whenever $\Mm$ is finite. \Qed

Theorem~\ref{notrans-thm} below serves as a contrast to
Theorem~\ref{trans2-thm}.

\begin{thm}\label{notrans-thm} There is a finite binary model that is not 3-equivalent
to any transitive model.
\end{thm}

\noindent {\bf Proof.}
The binary model $\Mm$ has 45 elements and 4 basic relations
$S,G,R,B$. Let $\langle \mathbb{Z}_5,+\rangle$ denote the group of
non-negative numbers smaller than 5 with addition modulo 5, and let
$9=\{ 0,1,\dots,8\}$ denote the set of non-negative integers smaller
than 9. We define
\[ M = 5\times 9 .\]
\begin{description}
\item{} Let $s,g$ be permutations of 9 defined, in cycle form, as\\
$s=(012)(345)(678)$\quad and\quad $g=(136)(147)(258)$,
\item{} and let $r,b\subseteq 9\times 9$ be defined as\\
$r=\{0,3,6\}\times\{0,1,2\}\cup
\{1,4,7\}\times\{3,4,5\}\cup\{2,5,8\}\times\{6,7,8\}$ and\\
$b=\{0,4,8\}\times\{0,5,7\}\cup
\{1,5,6\}\times\{1,3,8\}\cup\{2,3,7\}\times\{2,4,6\}$.
\end{description}
Now, the basic relations of $\Mm$ are defined as
\begin{description}
\item{}
$S=\{\langle (i,j),(i,s(j))\rangle : i\in 5, j\in 9\}$,
\item{}
$G=\{\langle (i,j),(i,g(j))\rangle : i\in 5, j\in 9\}$,
\item{}
$R=\{\langle (i,j),(i+1,k)\rangle : i\in 5, (j,k)\in r\}$,
\item{}
$B=\{\langle (i,j),(i+2,k)\rangle : i\in 5, (j,k)\in b\}$.
\end{description}

We show that $\Mm$ is not 3-equivalent to any transitive model. Let
us call a model \emph{3,1-transitive} when any two elements of the
same 3-type can be taken to each other by an automorphism. A
transitive model $\Nn$ is 3,1-transitive, because let $a,b\in N$
have the same 3-types, then they have the same types in $\Nn$,
therefore there is an automorphism taking $a$ to $b$, by
transitivity of $\Nn$. Thus, it is enough to show that if $\Mm$ is
3-equivalent to $\Nn$ then $\Nn$ is not 3,1-transitive.

Assume that $\Mm$ is 3-equivalent to $\Nn$. For a first-order
formula $\rho(x,y)$ with free variables among $x,y$ let $\rho(\Mm)$
denote the relation that $\rho$ defines in $\Mm$, i.e.,
$\rho(\Mm)=\{(a,b) : \Mm\models\rho[a,b]\}$. Let $\RA{\Mm}$ be the
relation algebra of FO$3$-definable binary relations of $\Mm$, i.e.,
the universe of $\RA{\Mm}$ is $\{ \rho{(\Mm)} : \rho(x,y)\in
FO$3$\}$, and the operations of $\RA{\Mm}$ are the operations of
taking union, converse and relation composition of binary relations
together with the (base-sensitive) operations of taking complement
in $M\times M$ and the identity constant $\{ (u,u) : u\in M\}$ on
$M$.
Let $\RA{\Nn}$ denote the similar algebra of FO$3$-definable binary
relations of $\Nn$. We show the following:
\[\tag{1} \mbox{$\Mm$ is 3-equivalent to $\Nn$ implies that
$\RA{\Mm}$ is isomorphic to $\RA{\Nn}$.}\] Indeed, it is easy to
check that the relation $\{(\rho(\Mm),\rho(\Nn)) : \rho(x,y)\in
FO$3$\}$ is an isomorphism between $\RA{\Mm}$ and $\RA{\Nn}$ when
$\Mm$ is 3-equivalent to $\Nn$.

A \emph{base-automorphism} of $\RA{\Nn}$ is a permutation $\alpha$
of $N$ that leaves all elements of $\RA{\Nn}$ fixed when taking $Z$
to $\{ (\alpha(u),\alpha(v)) : (u,v)\in Z\}$. Now, $\RA{\Nn}$ is
called \emph{c-permutational} iff any element of $N$ can be taken to
any other by a base-automorphism.
\[\tag{2} \mbox{$\Nn$ is 3,1-transitive implies that $\RA{\Nn}$ is
c-permutational}.\]
To check (2), notice first that all elements in $\Nn$ have the same
3-type. This is so because $\Mm$ is such and this property can be
expressed with FO$3$ formulas $\{\exists x\rho(x,x)\to\forall
x\rho(x,x) : \rho(x,y)\in FO3\}$. Therefore, $\Nn$ is 3,1-transitive
means that each element of $N$ can be taken to any other element of
$N$ by an automorphism of $\Nn$. Finally, an automorphism $\alpha$
of $\Nn$ is a base-automorphism of $\RA{\Nn}$ and we are done.

From now on, we will use \cite{ADN}.
\[\tag{3} \mbox{$\RA{\Mm}$ is the algebra $\Aa$ defined in
\cite[section 2]{ADN}.}\]
Indeed, it can be checked that the basic relations $S,G,R,B$ of
$\Mm$ coincide with the relations $s,g,r_0,b_0$ in \cite[section
2]{ADN}. It is stated in \cite[p.375, line 16]{ADN} that $\Aa$ is
generated by these four elements, so each element of $A$ is
FO$3$-definable in $\Mm$. In the other direction, it is a theorem of
relation algebra theory that all FO$3$-definable elements of $\Mm$
can be generated from the basic relations of $\Mm$ with the
operations of $\RA{\Mm}$, see e.g., \cite[sec.3.9]{TG} or
\cite[Thm.3.32]{HH}. Thus, the elements of $A$ are exactly the
FO$3$-definable binary relations of $\Mm$ and we are done.

In the proof of \cite[Theorem 1]{ADN} it is proved that $\Aa$ is not
isomorphic to any c-permutational algebra, and so $\Nn$ cannot be
3,1-permutational by (1) and (2). The proof of
Theorem~\ref{notrans-thm} is complete. \Qed

\section{Two-variable fragment of FO has the weak Beth definability
property}\label{wbeth-sec}

We recall the definition of when the two-variable fragment FO$2$ has
the weak Beth definability property (wBDP).

\begin{Def}[wBDP for FO$2$]\label{wbeth-def}
Let $\Ll$ be a language with relation symbols of rank 2, let $\Th$
be any set of formulas of FO$2$($\Ll$), the set of formulas of
language $\Ll$ that contain, bound or free, only the variables
$x,y$. Assume that $\R$ is a binary relation symbol not occurring in
$\Ll$, let $\Ll^+$ denote $\Ll$ expanded with $\R$. Let $\Sigma(\R)$
be a set of formulas of FO$2$($\Ll^+$).
\begin{description}
\item{(i)}
We say that $\Sigma(\R)$ is a \emph{strong implicit definition of
$\R$ w.r.t.\ $\Th$} when in each model $\Mm$ of $\Th$ there is
exactly one relation $R\subseteq M\times M$ such that
$\langle\Mm,R\rangle\models\Sigma(\R)$. We say that $\Sigma(\R)$ is
just a \emph{weak implicit definition of $\R$ w.r.t.\ $\Th$} when in
each model $\Mm$ of $\Th$ there is at most one relation $R$ such
that $\langle\Mm,R\rangle\models\Sigma(\R)$. That is, with weak
definitions it is allowed that in some models there is no relation
at all satisfying the implicit definition.
\item{(ii)}
We say that $\Sigma(\R)$ \emph{can be made explicit w.r.t.\ $\Th$},
or that \emph{$\R$ has an explicit definition over $\Th$} when there
is a formula $\varphi\in FO2(\Ll)$ such that
$\Th\cup\Sigma(\R)\models \forall
x,y(\R(x,y)\leftrightarrow\varphi)$. In this case, we say that
$\varphi$ is an explicit definition of $\R$ in $\Th$.
\item{(iii)}
Now, \emph{FO$2$ has the weak Beth definability property} means that
each strong implicit definition of FO$2$ can be made explicit. \Qed
\end{description}
\end{Def}
It is proved in \cite{ACMNS} that there is a weak implicit
definition in FO$2$ that cannot be made explicit. The question is
whether there is also a strong implicit definition that cannot be
made explicit. The following theorem gives a negative answer to
this.

\begin{thm}\label{wbeth-thm}
FO$2$ has the weak Beth definability property.
\end{thm}

The proof of Theorem~\ref{wbeth-thm} uses the following two lemmas.
We say that $R\subseteq M\times M$ \emph{does not cut 2-types in
$\Mm$} if for all $a,b,c,d\in M$ we have ($R(a,b)$ iff $R(c,d)$)
whenever the 2-type of $(a,b)$ is the same as that of $(c,d)$ in
$\Mm$.

\begin{lem}[No-cut lemma]\label{nocut-lem} Assume that $\Sigma$ is a weak
implicit definition of $\R$ in $\Th$ and $\Th\cup\Sigma$ consists of
FO$2$-formulas. Assume that $\Mm$  is transitive and
$\langle\Mm,R\rangle\models\Th\cup\Sigma$. Then $R$ does not cut
2-types in $\Mm$.
\end{lem}

\noindent{\bf Proof of Lemma~\ref{nocut-lem}.} Let $\Th, \Sigma,
\Mm, R$ be as in the first two sentences of the statement of
Lemma~\ref{nocut-lem}. Assume that $R$ cuts a 2-type in $\Mm$. By
using this, we are going to define a relation $S$ distinct from $R$
which also satisfies $\Sigma$ in $\Mm$, contradicting that $\Sigma$
is an implicit definition.

Assume that \[\tag{a} \mbox{$R(a,b)$ and not $R(c,d)$ in $\Mm$}\]
for some $a,b,c,d$ such that $\Type(a,b,\Mm)=\Type(c,d,\Mm)$. Let
$T=\Type(a,b,\Mm)$ and
\[ t = \{ (m,n)\in M\times M : \Type(m,n,\Mm)=T\} .\]
We note that
\[\tag{d} \mbox{$T$ is not an identity type}\]
by (a) and $\Mm$ being transitive:  assume $(m,m)\in t$ for some
$m\in M$, then $(a,a), (c,c)\in t$ and there is an automorphism
$\alpha$ of $\Mm$ taking $a$ to $c$. Hence $(a,a)\in R$ iff
$(c,c)\in R$ because automorphisms preserve meanings of formulas and
hence they leave solutions of implicit definitions fixed. However,
$a=b$ and $c=d$ by $t$ being an identity type, contradicting (a).

For a binary relation $Z$, let $Z^{-1}=\{ (v,u) : (u,v)\in Z\}$
denote its inverse, $Z$ is symmetric means that $Z=Z^{-1}$. We
define $S\subseteq M\times M$ as follows. If $T\ne\inv{T}$ or $t\cap
R$ is symmetric, then we define $S$ by ``interchanging" $t\cap R$
with $t\setminus R$, i.e.,
\[ S=(R\setminus t)\cup (t\setminus R) .\]
If $T=\inv{T}$ and $t\cap R$ is not symmetric then $S$ is defined by
``interchanging" $(t\cap R\setminus R^{-1})$ with its inverse
$(t\cap R^{-1}\setminus R)$, i.e.,
\[ S=(R\setminus [t\cap R\setminus R^{-1}])\cup (t\cap R^{-1}\setminus R) .\]
Then $S$ is distinct from $R$ by $t\cap R$ being nonempty in the
first case, and by $(t\cap R\setminus R^{-1})\cup(t\cap
R^{-1}\setminus R)$ being nonempty in the second case. However,
\[\tag{s} \mbox{$S$ and $R$ differ only inside $t$},\]
that is, $[R(m,n)$ iff $S(m,n)]$ for all $(m,n)\in M\times
M\setminus t$.  We are going to show that
$\langle\Mm,S\rangle\models\Sigma$.
We define $J\subseteq (M\times M)\cup (M^2\times M^2)$ by requiring
for all $m,n,p,q\in M$ that
\begin{description}
\item{}
$\langle m,p\rangle\in J\quad\mbox{ iff }\quad
\Type(m,m,\Mm)=\Type(p,p,\Mm)$,
\item{}
$\langle (m,n),(p,q)\rangle\in J$\quad\mbox{ iff }\quad\\
$\Type(m,n,\Mm)=\Type(p,q,\Mm)$\mbox{ and }\\
$[(m,n)\in R\leftrightarrow (p,q)\in S)]$\mbox{ and
}\\
$[(n,m)\in R\leftrightarrow (q,p)\in S)]$.
\end{description}
We now show that $J$ is a 2-isomorphism between
$\langle\Mm,R\rangle$ and $\langle\Mm,S\rangle$. We check properties
(i)-(iv) in the definition of a 2-isomorphism. (i) is satisfied by
the definition of $J$. Restriction property (iii) is satisfied,
because if $\Type(a,b,\Mm)=\Type(p,q,\Mm)$ then
$\Type(a,a,\Mm)=\Type(b,b,\Mm)$ and $\Type(p,p,\Mm)=\Type(q,q,\Mm)$.

Checking local isomorphism property (ii): Assume
$\langle(m,n),(p,q)\rangle\in J$. Then $J$ is a local isomorphism
with respect to the language of $\Mm$ by \newline
$\Type(m,n,\Mm)=\Type(p,q,\Mm)$. We now check local isomorphism with
respect to the new relation symbol $\R$. Indeed, $R(m,m)$ iff
$S(p,p)$ holds because $R(m,m)$ iff $R(p,p)$ by the restriction
property, and $R(p,p)$ iff $S(p,p)$ by (s) and (d).  We have
$R(m,n)$ iff $S(p,q)$ and $R(n,m)$ iff $S(q,p)$  by the definition
of $J$. Thus property (ii) holds.

Checking back-and-forth property (iv): The first line is satisfied
by $\langle a,a\rangle\in J$ for all $a\in M$. For checking the
second line, let $m,n,m'\in M$, $\langle m,n\rangle\in J$. Let
$\alpha$ be an automorphism of $\Mm$ that takes $m$ to $n$. There is
such an $\alpha$ because $\langle m,n\rangle\in J$ and $\Mm$ is
transitive.

Assume that $\Type(m,m',\Mm)\notin\{ T,\inv{T}\}$. Let
$n'=\alpha(m')$. Then
\newline $\Type(n,n',\Mm)=\Type(m,m',\Mm)\notin\{ T,\inv{T}\}$,
because automorphisms do not change 2-types. By (s) then $R(e,f)$
iff $S(e,f)$ for all $(e,f)$ in $\{(m,m'),(m',m),\\
(n,n'),(n',n)\}$, thus $\langle(m,m'),(n,n')\rangle\in J$ by the
definition of $J$.

Assume now that $\Type(m,m',\Mm)=T$ and $T\ne\inv{T}$ or $t\cap R$
is symmetric. If $R(m,m')$, then let $\beta$ be an automorphism that
takes $c$ to $n$ and let $n'=\beta(d)$. Then $(n,n')\in (t\setminus
R)$ by (a), so $(n,n')\in S$ by the definition of $S$. Assume
$T\ne\inv{T}$, then $(m',m)\in R$ iff $(n',n)\in S$ by (s). Assume
that $T=\inv{T}$ and $t\cap R$ is symmetric. Then $t\cap R^{-1}$ is
also symmetric, hence $(m',m)\in R$ and $(n',n)\notin R$, so
$(n',n)\in S$. Thus $\langle (m,m'),(n,n')\rangle\in J$. If not
$R(m,m')$, then let $\beta$ be an automorphism that takes $a$ to $n$
and let $n'=\beta(b)$. From here on, the argument showing
$\langle(m,m'),(n,n')\rangle\in J$ is completely analogous to the
previous case.

Assume now that $\Type(m,m',\Mm)=T=\inv{T}$ and $t\cap R$ is not
symmetric, say, $(e,f)\in t\cap R\setminus R^{-1}$. If $(m,m')\in
(R\cap R^{-1})$ then let $n'=\alpha(m')$. Then $(n,n')\in (R\cap
R^{-1})$ by $\alpha(R)=R$ and so $\alpha(R^{-1})=R^{-1}$. Also,
$(n,n')\in S\cap S^{-1}$ in this case, by the definition of $S$.
Hence $\langle (m,m'),(n,n')\rangle \in J$. The case $(m,m')\notin
(R\cup R^{-1})$ is completely analogous. Assume $(m,m')\in
R\setminus R^{-1}$. Let $\beta$ be an automorphism taking $f$ to
$n$, and let $n'=\beta(e)$. Then $(n,n')\in R^{-1}\setminus R$, so
$(n,n')\in S\setminus S^{-1}$ by the definition of $S$, and so
$\langle (m,m'),(n,n')\rangle\in J$ by the definition of $J$. The
case $(m,m')\in R^{-1}\setminus R$ is completely analogous.

The case when $\Type(m,m',\Mm)=\inv{T}$ can be proved in a
completely analogous way. By this, checking the second line is
finished.
The third line of (iv) follows in our case from the second one by
noticing that $J$ is symmetric both in $M\times M$ and in $M^2\times
M^2$.

Thus $J$ is a 2-isomorphism between $\langle\Mm,R\rangle$ and
$\langle\Mm,S\rangle$, so $\langle\Mm,S\rangle\models\Sigma$. Thus
both $R$ and the distinct $S$ satisfy $\Sigma$, this contradicts
$\Sigma$ being an implicit definition. The proof of
Lemma~\ref{nocut-lem} is complete. \Qed

The next lemma is a kind of characterization of those implicit
definitions that cannot be made explicit.

\begin{lem}[Cut lemma]\label{cut-lem}
Assume that $\Sigma$ is a weak implicit definition of $\R$ in $\Th$,
and $\Th\cup\Sigma$ consists of FO$2$ formulas. Statements (i) -
(iii) below are equivalent.
\begin{description}
\item[(i)] $\Sigma$ can be made explicit in $\Th$ by a
FO$2$-formula.
\item[(ii)] $R$ does not cut 2-types in $\Mm$
whenever $\langle\Mm,R\rangle\models\Th\cup\Sigma$.
\item[(iii)] $R$ does not cut 2-types in $\Mm$ whenever
$\langle\Mm,R\rangle\models\Th\cup\Sigma$ and $\Mm$ is
2-homogeneous.
\end{description}
\end{lem}

\noindent {\bf Proof of Lemma~\ref{cut-lem}.} Clearly, (i) implies
(ii) and (ii) implies (iii).
To show that (iii) implies (i), assume that (i) does not hold, we
will infer the negation of (iii).

We will refer to the negation of (i) just as ``$\R$ is not
2-definable". Thus we want to show that there is a model
$\langle\Mm,R\rangle\models\Th\cup\Sigma$ such that $\Mm$ is
2-homogeneous and $R$ cuts a 2-type in $\Mm$.

By a \emph{2-partition} we understand a system $\langle\pi_i : i\le
n\rangle$ of FO$2$-formulas in the language of $\Th$ such that
$\Th\cup\Sigma\models (\bigvee\{ \pi_i : i\le n\}\land\bigwedge\{
\lnot(\pi_i\land\pi_j) : i<j\le n\}$. We say that \emph{$\R$ cannot
cut $\pi$} when $\Th\cup\Sigma\models \forall x,y(\pi(x,y)\to
\R(x,y))\lor(\forall x,y(\pi(x,y)\to\lnot \R(x,y))$. We say that
\emph{$\R$ cannot cut into the 2-partition} $\langle\pi_i : i\le
n\rangle$ when $\R$ cannot cut any $\pi_i$ for $i\le n$.
We will use the following statement
\[ \tag{L1}\label{L1} \mbox{$\R$ is 2-definable\quad iff\quad there is a 2-partition $\R$ cannot cut into}.\]
Indeed, if $\R$ is definable by the FO$2$-formula $\rho(x,y)$, then
$\R$ cannot cut into $\langle \rho, \lnot\rho\rangle$. In the other
direction, assume that $\R$ cannot cut into $\langle\pi_i : i\le
n\rangle$. Let us treat natural numbers in the von Neumann's sense,
i.e., each natural number $n$ is the set of smaller natural numbers:
$n=\{ i\in\omega : i<n\}$. For $J\subseteq n$ let
$\pi(J)=\bigwedge\{ \pi_j : j\in J\}\land\bigwedge\{\lnot\pi_j :
j\in n\setminus J\}$. By the assumption that $\R$ cannot cut into
$\langle\pi_i : i<n\rangle$ we have that $\R$ is a union of some
$\pi_i$s in each model of $\Th\cup\Sigma$, i.e.,
$\Th\cup\Sigma\models \bigvee\{ \R(x,y)\leftrightarrow\pi(J) :
J\subseteq n\}$. Since $\Sigma$ is an implicit definition, we have
$\Th\cup\Sigma(\R)\cup\Sigma(\R')\models \R(x,y)\leftrightarrow
\R'(x,y)$ where $\R'$ is a brand new binary relation symbol, so by
compactness
\[ \tag{s}\label{s}\Th\cup\Sigma_0(\R)\cup\Sigma_0(\R')\models \R(x,y)\leftrightarrow
\R'(x,y)\quad\mbox{for some finite $\Sigma_0\subseteq\Sigma$.}\] Let
$\sigma(\R)=\bigwedge\Sigma_0(\R)$ for a $\Sigma_0$ satisfying (s).
Assume that $\langle\Mm,R\rangle\models\Th\cup\Sigma$. Then there is
a unique $J$ such that $\Mm\models \R(x,y)\leftrightarrow\pi(J)$,
since $\pi(J)\land\pi(K)$ is inconsistent for distinct $J$ and $K$.
By $\langle\Mm,R\rangle\models\Sigma$ and $\Sigma_0\subseteq\Sigma$
we have that $\Mm\models\sigma(\pi(J))$. Also,
$\Mm\models\sigma(\pi(K))$ is not true for $K\ne J$ by (s) since
$\pi(J)$ and $\pi(K)$ define distinct relations in $\Mm$. This shows
that $\Th\cup\Sigma\models
\R(x,y)\leftrightarrow\bigwedge\{\sigma(\pi(J))\to\pi(J) :
J\subseteq n\}$, and this is a 2-definition for $\R$. Statement (L1)
has been proved.

To prove (iii), we construct a 2-type $T$ in the language of $\Th$
such that the set
\[ \tag{L2}\label{L2} \Th\cup\Sigma\cup \{ \R(x,y), \lnot \R(z,v)\}\cup T(x,y)\cup
T(z,v)\] of FO-formulas is consistent. That $T$ is a 2-type means
that either $\rho\in T$ or $\lnot\rho\in T$ for all FO$2$-formulas
$\rho$ in the language of $\Th$. We say that a set of open
FO-formulas is consistent when there is a model and an evaluation
that make the set true. Equivalently, one can consider the free
variables in the set to be constants, we will use this second
option.

Let $\tau$ be a FO$2$-formula in the language of $\Th$. A
\emph{2-partition of $\tau$} is a system $\langle\pi_i : i\le
n\rangle$ such that
$\Th\cup\Sigma\models(\tau\leftrightarrow\bigvee\{\pi_i :
i<n\})\land\bigwedge\{\lnot(\pi_i\land\pi_j) : i<j<n\}$.
Let $T$ be a set of FO$2$ formulas in the language of $\Th$. We say
that \emph{$T$ is good} iff $\R$ can cut into any 2-partition of
$\bigwedge T_0$, for all finite subsets $T_0$ of $T$.
Clearly, a directed union of good sets is a good set again, so there
is a maximal one among the good sets by Zorn's lemma. We will show
that any maximal good set is a 2-type and that the set in (L2) with
any good $T$ is consistent. We begin with this second statement.

Assume that $T$ is good and let $T_0\subseteq T$ be finite. Then
$\R$ can cut into any 2-partition of $\bigwedge T_0$, in particular
$\R$ can cut into $\bigwedge T_0$. This means that there is a model
of $\Th\cup\Sigma\cup\{ \R(x,y), \lnot \R(z,v)\}\cup T_0(x,y)\cup
T_0(z,v)$. By compactness, the set in (L2) is consistent.

To show that a maximal good $T$ is a 2-type, let $T$ be any good set
and let $\rho$ be any FO$2$ formula in the language of $\Th$. We
show that either $T\cup\{\rho\}$ is good or $T\cup\{\lnot\rho\}$ is
good. Assume that neither of $T\cup\{\rho\}$ and
$T\cup\{\lnot\rho\}$ is good. Then there are finite subsets $T_0,
T_1$ of $T$ and 2-partitions $\pi=\langle \pi_i : i\le n\rangle$ of
$\rho\land\bigwedge T_0$ and $\delta=\langle \delta_j : j<m\rangle$
of $\lnot\rho\land\bigwedge T_1$ such that $\R$ cannot cut into
either of these two partitions. We can now combine $\pi$ and
$\delta$ to form a 2-partition $\sigma$ of $\bigwedge (T_0\cup T_1)$
by letting the members of the partition $\sigma$ be
$\pi_i\land\bigwedge T_1$ and $\delta_j\land\bigwedge T_0$ for $i<n,
j<m$. Clearly, $\R$ cannot cut into $\sigma$ by our assumption that
$\R$ cannot cut into either of $\pi$ and $\delta$; this contradicts
to $T$ being good. With this, we have proved that any maximal good
$T$ is a 2-type.

By the above, we now have a 2-type $T$ such that the set
$\Delta=\Th\cup\Sigma\cup\{ \R(x,y), \lnot \R(z,v)\}\cup T(x,y)\cup
T(z,v)$ is consistent. Let then $\langle \Mm,R\rangle$ be any
$\omega$-saturated model of $\Delta$. Then $\Mm$ is also
$\omega$-saturated, and so it is 2-homogeneous by statement (S) in
the proof of Theorem~\ref{trans2-thm}. Also, $\R$ cuts the 2-type
$T$ in $\Mm$ by $\langle\Mm,R\rangle\models\Delta$. We derived the
negation of (iii) from the negation of (i), and this finishes the
proof of Lemma~\ref{cut-lem}. \Qed

\noindent{\bf Proof of Theorem~\ref{wbeth-thm}.}
Assume that $\Sigma$ is a strong implicit definition of $\R$ w.r.t.\
$\Th$. We are going to show that $\Sigma$ can be made explicit,
i.e., $\R$ has an explicit definition over $\Th$ that uses only two
variables.

Take any 2-homogeneous model $\Mm$ of $\Th$, and let
$\overline{\Mm}$ be a transitive model with $I$ a 2-isomorphism
between $\Mm$ and $\overline{\Mm}$. There are such $\overline{\Mm}$
and $I$ by Theorem~\ref{trans-thm}. $\overline{\Mm}$ is a model of
$\Th$ because 2-isomorphic models satisfy the same FO$2$-formulas.
Since $\Sigma$ is a \emph{strong} implicit definition of $\R$ in
$\Th$, there is $\overline{R}$ which satisfies $\Sigma$ in
$\overline{\Mm}$, i.e.,
$\langle\overline{\Mm},\overline{R}\rangle\models\Sigma$.

Since $\overline{\Mm}$ is transitive, $\overline{R}$ does not cut
2-types in $\overline{\Mm}$, by Lemma~\ref{nocut-lem}. We now
``transfer" $\overline{R}$ to the model $\Mm$ by the following
definition: take any pair $(a,b)$ in $\Mm$. We define $R$ so that
this pair is related by $R$ if and only if there is a pair of the
same type in $\overline{\Mm}$ which is related by $\overline{R}$.
Formally:
\[ R = \{(a,b)\in M\times M : \exists c,d[\overline{R}(c,d)\mbox{
and }\Type(a,b,\Mm)=\Type(c,d,\overline{\Mm})\} .\]

We now show that with this definition, the 2-isomorphism $I$ between
$\Mm$ and $\overline{\Mm}$ remains a 2-isomorphism between
$\langle\Mm,R\rangle$ and
$\langle\overline{\Mm},\overline{R}\rangle$. Since $I$ is a
2-isomorphism between $\Mm$ and $\overline{\Mm}$, it satisfies
conditions (i), (iii), (iv) in the definition of a 2-isomorphism,
and it also satisfies condition (ii) for atomic formulas other than
$R(v,z)$. Therefore, we only have to show that if $\langle (a,b),
(a',b')\rangle\in I$, then $R(a,b)$ iff $\overline{R}(a',b')$.

Assume that $\langle (a,b), (a',b')\rangle\in I$. If $R(a,b)$, then
there are $c,d\in\overline{\Mm}$ such that $\overline{R}(c,d)$ and
$\Type(c,d,\overline{\Mm})=\Type(a,b,\Mm)$, by the definition of
$R$.
The 2-type of $(a',b')$ is also the same as that of $(a,b)$, since
they are $I$-related by assumption. Hence the 2-type of $(c,d)$ is
the same as that of $(a',b')$ (since they both equal the 2-type of
$(a,b)$). By $\overline{R}(c,d)$ we now  get $\overline{R}(a',b')$,
since we have seen that $\overline{R}$ does not distinguish elements
of the same 2-type. In the other direction, assume that
$\overline{R}(a',b')$. Since $(a,b)$ is $I$-related to $(a',b')$,
their 2-types equal, hence $R(a,b)$ by the definition of $R$. By
this, we have seen that $I$ is a 2-isomorphism between the expanded
models $\langle\Mm,R\rangle$ and
$\langle\overline{\Mm},\overline{R}\rangle$. Since the latter is a
model of $\Sigma$, we get that $\langle\Mm,R\rangle$ is a model of
$\Sigma$, too. By its definition, $R$ does not cut 2-types in $\Mm$.

Since $\Sigma$ is a weak definition over $\Th$ and the 2-homogeneous
$\Mm\models\Th$ was chosen arbitrarily, we get that $R$ does not cut
2-types in $\Mm$ whenever $\langle\Mm,R\rangle\models\Th\cup\Sigma$
and $\Mm$ is 2-homogeneous. Thus, $\Sigma$ can be made explicit in
$\Th$, by Lemma~\ref{cut-lem} and this finishes the proof of
Theorem~\ref{wbeth-thm}. \Qed

\bigskip
\noindent{\bf On some connections with algebra and the literature.}

(1) There is an underlying algebraic intuition behind the proofs in
the paper. For example, one may wonder about the role of
2-homogeneity that shows up in Theorem~\ref{trans-thm}. For any
model of any logic, the concepts (that is, explicitly defined
notions) form a natural algebra with the connectives of the logic as
operations, this is called the \emph{concept algebra} of the model.
The ``types" of a model form a similar natural algebra only when
some kind of homogeneity is satisfied. In the case of FO$2$, this
condition is 2-homogeneity. The algebra of types is an atomic
superalgebra of the concept algebra.
A 2-partial isomorphism between binary models induces an isomorphism
between the respective algebras of types, and vice versa, any
isomorphism between algebras of types induces a 2-isomorphism
between the respective models. All this is part of a general
approach to the ``algebra behind logic", for more on this see
\cite[Part II]{ANSHbPhL01}.

(2) A bridge between logic and algebra is elaborated in, e.g.,
\cite{BP89}, \cite[section 4.3]{HMTII}, or \cite[Part
II]{ANSHbPhL01}. In this bridge, a class Alg(L) of algebras is
associated to any (decent) logic L. Namely, Alg(L) is the infinitary
quasi-equational hull of the class CA(L) of all concept algebras of
models of logic.
This correspondence is used for stating equivalence theorems of the
kind: L has logical property $LP$ if and only if Alg(L) has
algebraic property $AP$. A satisfying fact is that, usually, to
natural logical properties natural algebraic properties correspond
this way (see the first sentence of \cite{P72}).

(3) To Beth definability BDP of a logic the corresponding algebraic
property is surjectivity of epimorphisms in Alg(L) considered as a
category (ES), see \cite[Thm.5.6]{HMTII},
\cite[Thm.6.11(i)]{ANSHbPhL01}. Indeed, failure of the BDP for
$n$-variable logics was proved first via showing that ES fails in
the category of their concept algebras, see \cite{ACMNS}. Failure of
BDP for equality-free 2-variable logic is proved also in algebraic
form in \cite{SBMP}.

(4) Weak Beth definability wBDP possesses several natural algebraic
equivalent properties, see \cite[Thm.6.11(iii) and
Thm.6.12]{ANSHbPhL01}  and \cite{Hoo00}. Three of the corresponding
algebraic properties are (i) ``surjectivity of CA(L)-extendible
epimorphisms", (ii) ``Alg(L) is the smallest full reflective
subcategory containing MA(L)", and (iii) ``MA(L) generates Alg(L) by
taking limits of diagrams of algebras". Here, MA(L) is the class of
maximal (with respect to containment) members of CA(L).

(5) Craig Interpolation Property CIP and Beth definability property
BDP are related, for example, BDP is often proved from CIP. The
interpolation properties also have several variants that coincide in
the case of FO while happen to be distinct in other logics, e.g., in
some modal logic.
In the spirit of the bridge theorems, to some interpolation property
of a logic L some kind of amalgamation property of Alg(L)
corresponds, see, e.g., \cite[Thm.6.15]{ANSHbPhL01}. A landmark
paper on amalgamation and related properties in algebraic versions
of FO is \cite{P72}.

(6) The algebraic equivalent of CIP is strong amalgamation property
SAP. The superamalgamation property SUPAP of Alg(L) is proved to
correspond to a stronger version sCIP of Craig interpolation
property in \cite{Maksi, Mada97}.
The question whether they coincide for variants of FO is asked in
\cite{P72}. This question is answered in \cite{SS06}, where a
variant of $n$-variable logic, for any finite $n\ge 3$, is
constructed that has CIP but not sCIP. This is an analogous result
to the one in the present paper concerning Craig Interpolation
Property in place of Beth definability property.
We mention that all of the other questions in the two diagrams in
\cite{P72} have been answered in the meantime, see
\cite{MadaSaye06}.

(7) There is an extended literature connecting logic and algebra. A
tiny sample is \cite{Comer, Gyenis, KS, KLSF, Say53}.
\bigskip

In computer science, $n$-variable logic usually is used in a
stronger version where it is endowed with infinite conjunctions and
disjunctions, see \cite{Hodk}. Let $L^n_{\infty,\omega}$ denote this
logic as in \cite{Hodk}. Failure of wBDP for $L^n_{\infty,\omega}$
in the case of $n\ge 3$ is proved in \cite{Hodk}, but it also
follows from \cite{ANwBeth14} (since there the proof is based on
finite counterexamples), failure of BDP for $L^2_{\infty,\omega}$
follows from the proof in \cite{ACMNS}.

\begin{qes} Does $L^2_{\infty,\omega}$ have weak Beth definability
property?
\end{qes}

\section*{Acknowledgements} We thank Zal\'an Gyenis and G\'abor
S\'agi for enjoyable (transitive) discussions on the subject. We
also thank the referee for the many useful suggestions.

\bigskip\bigskip\bigskip

\noindent Alfr\'ed R\'enyi Institute of Mathematics\\
Budapest, Re\'altanoda st.\ 13-15, H-1053 Hungary\\
andreka.hajnal@renyi.hu, nemeti.istvan@renyi.hu

\end{document}